\font\tta=cmr10 scaled\magstep1

\font\ttc=cmss12

\font\tte=cmbsy10

\font\ttf=cmr8

\font\tth=cmbx10 scaled\magstep2

\def\proj#1{\hbox{{\tta \hskip 2pt P\hskip -10pt I~$_#1$}}}

\def\proja#1{\hbox{{\ttf \hskip 2pt P\hskip -8pt I~$_{\ttf #1}$}}}

\def\cmplxb{\ \hbox{{\ttc \hskip -2pt C\hskip -5pt I \hskip -1pt}}}
\def\cx#1{\cmplxb^{#1}}

\def\coker{\mathop{\rm coker}\nolimits}
\def\Hom{\mathop{\rm Hom}\nolimits}
\def\Ext{\mathop{\rm Ext}\nolimits}

\def\lra{\longrightarrow}

\def\sigg{\mathop{\hbox{$\displaystyle\sum$}}\limits}

\def\carre{$\Box$}
\def\timex{\setbox250=\hbox{$\times$}\hskip 5pt\Box\hskip -0.75em
{\raise 1.5pt\vbox{\box250}}\hskip 5pt}
\def\paragra{{\tte \char120}}
\def\para{\paragra~\hskip -2pt}

\def\hfl#1#2{\smash{\mathop{\hbox to 12mm{\rightarrowfill}}
\limits^{\scriptstyle#1}_{\scriptstyle#2}}}

\def\q#1#2{{{#1}\over{#2}}}
\def\m#1{{\hbox{$#1$}}}

\def\ot{\otimes}
\def\D{{\cal D}}
\def\E{{\cal E}}
\def\F{{\cal F}}

\def\O{{\cal O}}
\def\U{{\cal U}}
\def\V{{\cal V}}

\def\I{{\cal I}}
\def\P{{\cal P}}
\def\T{{\cal T}}
\def\SS{{\cal S}}
\def\og{\leavevmode\raise.3ex\hbox{$\scriptscriptstyle\langle\!\langle$}}
\def\fg{\leavevmode\raise.3ex\hbox{$\scriptscriptstyle\,\rangle\!\rangle$}}
\def\dem{\noindent{\em D\'emonstration}. }

\documentstyle[twoside,leqno,12pt]{article}
\begin{document}
\newtheorem{xprop}{Proposition}[section]
\newtheorem{xlemm}[xprop]{Lemme}
\newtheorem{xtheo}[xprop]{Th\'eor\`eme}
\newtheorem{xcoro}[xprop]{Corollaire}
\newtheorem{quest}{Question}
\newtheorem{defin}{D\'efinition}

\markboth{Jean-Marc Dr\'ezet}{Fibr\'es prioritaires g\'en\'eriques}
\title{Fibr\'es prioritaires g\'en\'eriques instables sur le plan projectif}
\author{Jean-Marc Dr\'ezet}
\date{}
\maketitle

\def\refname{R\'ef\'erences}
{\ttf
\noindent  \ Universit\'e Paris 7, UMR 9994 du CNRS, 
Aile 45-55, 5\m{^e} \'etage

\vskip -0.2cm

\noindent \ 2, place Jussieu, F-75251 Paris Cedex 05, France

\vskip -0.2cm

\noindent \ e-mail : drezet@math.jussieu.fr
}

\bigskip

\bigskip

\bigskip

{\parindent 3cm
{\tth \hskip 3cm Sommaire}

\bigskip

1 - Introduction

\medskip

2 - Fibr\'es exceptionnels

\medskip

3 - Fibr\'es prioritaires g\'en\'eriques

}

\bigskip

\bigskip

\section{Introduction}

Les faisceaux prioritaires sur \proj{2} ont \'et\'e
introduits par A. Hirschowitz et Y. Laszlo dans \cite{hi_la}. Rappelons qu'un
faisceau coh\'erent $\E$ sur \proj{2} est dit {\em prioritaire} s'il est sans
torsion et si \ \m{\Ext^2(\E,\E(-1))=0}. Par exemple les faisceaux semi-stables 
au sens de Gieseker-Maruyama sont prioritaires. On s'int\'eresse ici \`a la
structure pr\'ecise du faisceau prioritaire g\'en\'erique de rang $r$ et
de classes de Chern \m{c_1}, \m{c_2} lorsqu'il n'existe pas de faisceau 
semi-stable de m\^emes rang et classes de Chern.

D'apr\`es \cite{hi_la}, le {\em champ}
des faisceaux prioritaires est lisse et irr\'eductible. Les conditions 
d'existence des faisceaux prioritaires sont les suivantes : posons
$$\mu \ = \ \q{c_1}{r}, \ \ \ \Delta \ = \ \q{1}{r}(c_2 -
\q{r-1}{2r}c_1^2),$$
(si $\E$ est un faisceau coh\'erent $\E$ sur \proj{2} de rang $r$ et de
classes de Chern \m{c_1}, \m{c_2}, on appelle $\mu=\mu(\E)$ la {\em pente} de
$\E$ et $\Delta=\Delta(\E)$ le {\em discriminant} de $\E$). 
Alors, si \ \m{-1\leq\mu\leq 0}, il existe un faisceau prioritaire de
pente $\mu$ et de discriminant $\Delta$ si et seulement si on a
$$\Delta \ \geq \ - \q{\mu(\mu+1)}{2}.$$
Les conditions d'existence des faisceaux semi-stables sur \proj{2} sont 
rappel\'ees ci-dessous. On peut voir qu'il existe beaucoup de triplets
\m{(r,c_1,c_2)} tels qu'il existe un faisceau prioritaire de rang $r$ et 
de classes de Chern \m{c_1}, \m{c_2} mais pas de faisceau semi-stable avec les
m\^emes invariants. 

Les conditions d'existence des faisceaux semi-stables sur \proj{2} (cf \cite{dr_lp}) s'expriment en fonction des seules variables
$\mu$ et $\Delta$. On montre qu'il existe une unique fonction $\delta(\mu)$
telle qu'on ait \ \m{\dim(M(r,c_1,c_2)) > 0} \ si et seulement si \
\m{\Delta\geq\delta(\mu)}. La fonction \m{\delta(\mu)} est d\'ecrite \`a l'aide 
des {\it fibr\'es exceptionnels}.

On dit qu'un faisceau coh\'erent $\E$ sur \proj{2} est {\it exceptionnel} si
$\E$ est {\it simple} (c'est-\`a-dire si les seuls endomorphismes de $\E$ sont
les homoth\'eties), et si
$$\Ext^1(\E,\E) \ = \ \Ext^2(\E,\E) \ = \ \lbrace 0\rbrace.$$
Un tel faisceau est alors localement libre et stable, et la vari\'et\'e de
modules de faisceaux semi-stables correspondante contient l'unique point $\E$.
Il existe une infinit\'e d\'enombrable de fibr\'es exceptionnels, et un
proc\'ed\'e simple permet de les obtenir tous \`a partir des fibr\'es en
droites (cf. \cite{dr1}). Notons qu'un fibr\'e exceptionnel est 
uniquement d\'etermin\'e par sa pente.
Soit $F$ un fibr\'e exceptionnel. On note \m{x_F} la
plus petite solution de l'\'equation
$$X^2-3X+\q{1}{rg(F)^2} \ = \ 0.$$
Alors on montre que les intervalles \
\m{\rbrack\mu(F)-x_F,\mu(F)+x_F\lbrack} \  
constituent une partition de l'ensemble des nombres rationnels. On va d\'ecrire
la fonction \m{\delta(\mu)} sur cet intervalle. Posons
$$P(X) = \q{X^2}{2}+\q{3}{2}X+1.$$
Sur l'intervalle \ \m{\rbrack\mu(F)-x_F,\mu(F)\rbrack}, on a
$$\delta(\mu) \ = \ P(\mu-\mu(F))-\q{1}{2}(1-\q{1}{rg(F)^2}),$$
et sur \  \m{\lbrack\mu(F),\mu(F)+x_F\lbrack}, on a
$$\delta(\mu) \ = \ P(\mu(F)-\mu)-\q{1}{2}(1-\q{1}{rg(F)^2}).$$
On obtient les courbes $D(F)$ et $G(F)$ repr\'esent\'ees sur la figure qui suit.
Ce sont des segments de coniques. 

On consid\`ere maintenant la courbe \ \m{\Delta=\delta'(\mu)} \ d\'efinie de
la fa\c con suivante : sur l'intervalle \ 
\m{\rbrack\mu(F)-x_F,\mu(F)+x_F\lbrack}, on a
$$\delta'(\mu) = \delta(\mu) - \q{1}{rg(F)^2}(1-\q{1}{x_F}\mid\mu(F)-\mu\mid).$$
On obtient ainsi les segments de coniques $D'(F)$ et $G'(F)$. Le point
\m{(\mu(F),\delta'(\mu(F)))} est la paire \m{(\mu,\Delta)} correspondant au
fibr\'e exceptionnel $F$. Le point \m{(\mu(F),\delta(\mu(F)))} est le
sym\'etrique de $F$ par rapport \`a la droite \ \m{\Delta=1/2}. Notons que si
$\mu$ est un nombre rationnel diff\'erent de la pente d'un fibr\'e 
exceptionnel, le nombre $\delta'(\mu)$ est irrationnel.
Ces courbes,
sur l'intervalle \ \m{\rbrack\mu(F)-x_F,\mu(F)+x_F\lbrack} \ , sont
repr\'esent\'ees ci-dessous :

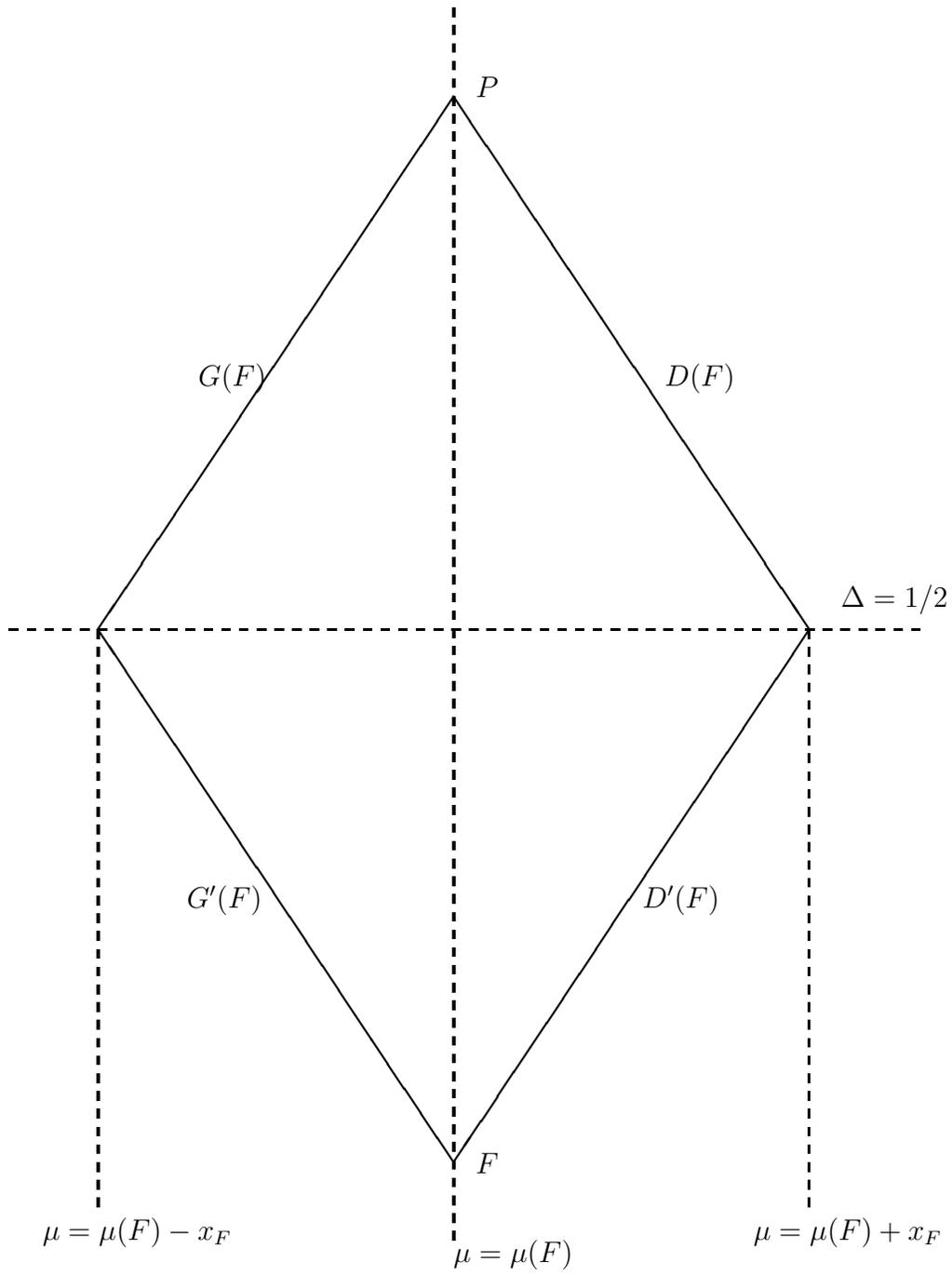
\vfill\eject

\setlength{\unitlength}{0.012500in}%
\begin{picture}(410,565)(200,235)
\thicklines
\multiput(400,800)(0.00000,-7.98561){70}{\line( 0,-1){  3.993}}
\multiput(610,520)(-7.96117,0.00000){52}{\line(-1, 0){  3.981}}
\put(400,760){\line(-2,-3){160}}
\put(400,760){\line( 2,-3){160}}
\put(560,520){\line(-2,-3){160}}
\put(400,280){\line(-2, 3){160}}
\multiput(240,520)(0.00000,-8.00000){33}{\line( 0,-1){  4.000}}
\multiput(560,520)(0.00000,-8.00000){33}{\line( 0,-1){  4.000}}
\put(285,630){\makebox(0,0)[lb]{\smash{$G(F)$}}}
\put(495,630){\makebox(0,0)[lb]{\smash{$D(F)$}}}
\put(280,395){\makebox(0,0)[lb]{\smash{$G'(F)$}}}
\put(485,395){\makebox(0,0)[lb]{\smash{$  D'(F)$}}}
\put(410,275){\makebox(0,0)[lb]{\smash{$F$}}}
\put(410,760){\makebox(0,0)[lb]{\smash{$P$}}}
\put(570,530){\makebox(0,0)[lb]{\smash{           $\Delta=1/2$}}}
\put(400,235){\makebox(0,0)[lb]{\smash{$\mu=\mu(F)$}}}
\put(215,245){\makebox(0,0)[lb]{\smash{$\mu=\mu(F)-x_F$}}}
\put(535,245){\makebox(0,0)[lb]{\smash{$\mu=\mu(F)+x_F$}}}
\end{picture}

\bigskip

\bigskip

\bigskip

Pour tout point $x$ de \proj{2}, soit \m{\I_x} le faisceau d'id\'eaux du point 
$x$. On a 
$$\Ext^1(\I_x,\O)\simeq\cx{}.$$ 
Soit \m{\V_x} l'unique faisceau
extension non triviale de \m{\I_x} par $\O$. On va d\'emontrer le

\vfill\eject

\noindent{\bf Th\'eor\`eme A : }{\em Soient $r$, \m{c_1}, \m{c_2} des entiers,
avec \m{r\geq 1}, \m{-1<\mu\leq 0}, 
$$\Delta \ \geq \ \q{\mu(\mu+1)}{2},$$
et tels que la vari\'et\'e \m{M(r,c_1,c_2)} soit vide.

\medskip

\noindent 1 - Si \ \m{\Delta < \delta'(\mu)}, il existe des fibr\'es
exceptionnels $E_0$, $E_1$, $E_2$, des espaces vectoriels de dimension finie
$M_0$, $M_1$, $M_2$, dont un au plus peut \^etre nul, tels que le faisceau
prioritaire g\'en\'erique de rang $r$ et de classes de Chern $c_1$, $c_2$ soit
isomorphe \`a 
$$(E_0\ot M_0)\oplus(E_1\ot M_1)\oplus(E_2\ot M_2).$$

\medskip

\noindent 2 - On suppose que \m{c_1\not = 0} ou \m{c_2>1}.
Si \ \m{\Delta > \delta'(\mu)}, soit $F$ l'unique fibr\'e 
exceptionnel tel que \ \m{\mu\in \ \rbrack\mu(F)-x_F,\mu(F)+x_F\lbrack}. Alors
si \ \m{\mu\leq\mu(F)}, l'entier
$$p \ = \ r.rg(F)(P(\mu-\mu(F))-\Delta-\Delta(F))$$
est strictement positif, et le faisceau prioritaire g\'en\'erique de rang $r$
et de classes de Chern $c_1$, $c_2$ est isomorphe \`a une somme directe
$$(F\ot \cx{p})\oplus\E,$$
o\`u $\E$ est un fibr\'e semi-stable situ\'e sur la courbe $G(F)$. De m\^eme,
si \ \m{\mu\geq\mu(F)}, l'entier
$$p \ = \ r.rg(F)(P(\mu(F)-\mu)-\Delta-\Delta(F))$$
est strictement positif, et le faisceau prioritaire g\'en\'erique de rang $r$
et de classes de Chern $c_1$, $c_2$ est isomorphe \`a une somme directe
$$(F\ot \cx{p})\oplus\E,$$
o\`u $\E$ est un fibr\'e semi-stable situ\'e sur la courbe $D(F)$.

\medskip

\noindent 3 - Si \ \m{c_1=0}, \m{c_2=1}, le faisceau prioritaire g\'en\'erique
de rang $r$ et de classes de Chern $c_1$, $c_2$ est isomorphe \`a une somme
directe du type
$$(\O\ot\cx{r-2})\oplus\V_x.$$}

\bigskip

Le r\'esultat pr\'ec\'edent apporte des pr\'ecisions sur ce
qui est d\'emontr\'e dans \cite{hi_la}, c'est-\`a-dire que s'il n'existe pas
des faisceau semi-stable de rang $r$ et de classes de Chern \m{c_1}, \m{c_2},
alors deux cas peuvent se produire : la filtration de Harder-Narasimhan du
faisceau prioritaire g\'en\'erique de rang $r$ et de classes de Chern \m{c_1},
\m{c_2} comporte deux termes, ou elle en comporte trois. Dans le premier cas,
un des termes est semi-exceptionnel (c'est-\`a-dire de la forme
\m{F\ot\cx{k}}, avec $F$ exceptionnel), et dans le second cas les trois termes
sont semi-exceptionnels.
 
\bigskip

Le th\'eor\`eme A permet aussi de conclure que s'il n'existe pas de faisceau
semi-stable de rang $r$ et de classes de Chern \m{c_1}, \m{c_2}, il n'existe
pas non plus d'{\em espaces de modules fins} de faisceaux de rang $r$ et de
classes de Chern \m{c_1}, \m{c_2} contenant au moins un faisceau prioritaire.

On appelle ici
{\it espace de modules fin} de faisceaux de rang $r$ et de classes de Chern
\m{c_1}, \m{c_2} sur \proj{2} la donn\'ee d'une vari\'et\'e alg\'ebrique lisse
$M$ non vide et d'un faisceau coh\'erent $\F$ sur \ \m{M\times\proj{2}}
poss\'edant les propri\'et\'es suivantes :

\medskip

\noindent (i) Le faisceau $\F$ est plat sur $M$ et pour tout point ferm\'e
$x$ de $M$, \ \m{\F_x=\F_{\mid \lbrace x\rbrace\times\proja{{\ttf 2}}}} est un
faisceau sans torsion sur \proj{2}, de rang $r$ et de classes de Chern
\m{c_1}, \m{c_2}.

\noindent (ii) Pour tout point ferm\'e $x$ de $M$, le faisceau $\F_x$ est
simple, on a \
\m{\Ext^2(\F_x,\F_x)=\lbrace 0\rbrace}, et le morphisme de d\'eformation
 infinit\'esimale de Koda\"ira-Spencer
$$T_xM\lra\Ext^1(\F_x,\F_x)$$
est surjectif.

\noindent (iii) Pour tous point ferm\'es distincts $x$ et $y$ de $M$, les
faisceaux \m{\F_x} et \m{\F_y} ne sont pas isomorphes.

\bigskip

Par exemple, si $r$, \m{c_1} et 
$$\chi \ = \ r - c_2 + \q{c_1(c_1+3)}{2}$$
sont premiers entre eux, et s'il existe un faisceau stable de rang $r$ et de
classes de Chern \m{c_1}, \m{c_2}, la vari\'et\'e de modules de ces faisceaux
stables, \'equip\'ee d'un {\it faisceau universel}, est un espace de
modules fin. Ceci sugg\`ere la conjecture suivante :

\bigskip 

\noindent{\bf Conjecture : }{\it Les seuls espaces de modules fins qui
soient projectifs sont les vari\'et\'es de modules de faisceaux stables,
lorsque $r$, \m{c_1} et $\chi$ sont premiers entre eux.}

\bigskip

Le th\'eor\`eme A entraine imm\'ediatement le 

\bigskip 

\noindent{\bf Th\'eor\`eme B : }{\em Soient $r$, \m{c_1}, \m{c_2} des entiers
avec \ $r\geq 1$. On suppose que la vari\'et\'e modules \m{M(r,c_1,c_2)} des
faisceaux semi-stables sur \proj{2} de rang $r$ et de classes de Chern \m{c_1},
\m{c_2} est vide. Alors il n'existe pas d'espace de modules fin de faisceaux
 de rang $r$ et de classes de Chern \m{c_1}, \m{c_2}, et contenant un faisceau
prioritaire.}

\bigskip

Il est possible de pr\'eciser le 1- du th\'eor\`eme A. On rappelle dans le
\para 2 la notion de {\em triade}, qui est un triplet particulier
\m{(E,F,G)} de fibr\'es exceptionnels. On ne consid\`ere ici que des
triades de fibr\'es exceptionnels dont les pentes sont comprises entre $-1$
et $0$. A la triade \m{(E,F,G)} correspond le {\em triangle} \m{\T_{(E,F,G)}}
du plan (de coordonn\'ees \m{(\mu,\Delta)}), dont les c\^ot\'es sont des
segments de paraboles et les sommets les points correspondant \`a $E$, $F$ et
$G$. Ce triangle est d\'efini par les in\'equations
$$\Delta\leq P(\mu-\mu(G))-\Delta(G), \ \
\Delta\geq P(\mu-\mu(H)+3)-\Delta(H), \ \
\Delta\leq P(\mu-\mu(E)+3)-\Delta(E),$$
$H$ \'etant le fibr\'e exceptionnel noyau du morphisme d'\'evaluation
\ \m{E\ot\Hom(E,F)\lra F}.
Soit {\bf T} l'ensemble des triades de fibr\'es exceptionnels dont les pentes
sont comprises entre $-1$ et $0$. Soit $\SS$ l'ensemble des points 
\m{(\mu,\Delta)} du plan tels que
$$-1\leq\mu\leq 0, \ \ -\q{\mu(\mu+1)}{2}\leq\Delta\leq\delta'(\mu).$$
On d\'emontrera le

\bigskip

\bigskip

\noindent{\bf Th\'eor\`eme C : }{\em 1 - Soient \m{(E,F,G)}, \m{(E',F',G')} des
\'el\'ements distincts de {\bf T}. Alors les triangles \m{\T_{(E,F,G)}} et
\m{\T_{(E',F',G')}} ont une intersection non vide si et seulement si cette
intersection est un sommet commun ou un c\^ot\'e commun. Dans le premier cas,
les fibr\'es exceptionnels correspondants sont identiques, et dans le second
les paires de fibr\'es exceptionnels correspondantes le sont.

\medskip

\noindent 2 - On a \ \ \ \m{\displaystyle
\SS\ = \ \bigcup_{(E,F,G)\in{\bf T}}\T_{(E,F,G)}}.

\medskip

\noindent 3 - Soient \m{r,c_1,c_2} des entiers, avec \ \m{r\geq 1},
$$\mu=\q{r}{c_1}, \ \ \Delta=\q{1}{r}(c_2-\q{r-1}{2r}c_1^2).$$
On suppose que \ \m{(\mu,\Delta)\in\T_{(E,F,G)}}. Soit $H$ le noyau du morphisme
d'\'evaluation \break \m{E\ot\Hom(E,F)\lra F}. Alors
$$m \ = \ r.rg(E).(P(\mu-\mu(E)+3)-\Delta(E)),$$
$$n \ = \ r.rg(H).(P(\mu-\mu(H)+3)-\Delta(H)),$$
$$p \ = \ r.rg(G).(P(\mu-\mu(G))-\Delta(G))$$
sont des entiers positifs ou nuls, et le fibr\'e prioritaire g\'en\'erique
de rang $r$ et de classes de Chern $c_1$, $c_2$ est de la forme
$$(E\ot\cx{m})\oplus(F\ot\cx{n})\oplus(G\ot\cx{p}).$$
}

\bigskip

\bigskip

\bigskip

\noindent{\bf Notations}

\medskip

Rappelons que le th\'eor\`eme de Riemann-Roch s'\'ecrit pour un faisceau
coh\'erent $E$ de rang positif sur \proj{2}
$$\chi(E) \ = \ rg(E).(P(\mu(E))-\Delta(E)),$$
\m{\chi(E)} d\'esignant la caract\'erisitique d'Euler-Poincar\'e de $E$.
Si $E$, $F$ sont des faisceaux coh\'erents sur \proj{2}, on pose
$$\chi(E,F) \ = \ \sigg_{0\leq i\leq 2}(-1)^i\dim(\Ext^i(E,F)).$$
On a, si \ \m{rg(E)>0} \ et \ \m{rg(F)>0},
$$\chi(E,F) \ = \ rg(E).rg(E).(P(\mu(F)-\mu(E))-\Delta(E)-\Delta(F)).$$
On a en g\'en\'eral, pour tout entier $i$, un isomorphisme canonique
$$\Ext^i(E,F) \ \simeq \ \Ext^{2-i}(F,E(-3))$$
(dualit\'e de Serre, cf. \cite{dr_lp}, prop. (1.2)).


\section{Fibr\'es exceptionnels}

\subsection{Construction des fibr\'es exceptionnels}

Les r\'esultats qui suivent ont \'et\'e d\'emontr\'es dans \cite{dr_lp} ou
\cite{dr1}. Un fibr\'e exceptionnel est enti\`erement d\'etermin\'e par sa
pente. Soit $\P$ l'ensemble des pentes de fibr\'es exceptionnels. Si 
\m{\alpha\in\P}, on note \m{E_\alpha} le fibr\'e exceptionnel de pente $\alpha$,
et \m{r_\alpha} son rang. On montre que \m{r_\alpha} et \m{c_1(E_\alpha)}
sont premiers entre eux. Soit \ \m{\Delta_\alpha=\Delta(E_\alpha)}. Alors on
a
$$\Delta_\alpha\ = \ \q{1}{2}(1-\q{1}{r_\alpha^2}),$$
(ce qui d\'ecoule du fait que \ \m{\chi(E_\alpha,E_\alpha)=1}). 

Soit $\D$ l'ensemble des nombres rationnels diadiques, c'est-\`a-dire pouvant
se mettre sous la forme \m{p/2^q}, $p$ et $q$ \'etant des entiers, \m{q\geq 0}.
On a une bijection
$$\epsilon : \E\lra\P.$$
Cette application est enti\`erement d\'etermin\'ee par les propri\'et\'es 
suivantes: 

\medskip

\noindent - Pour tout entier $k$, on a \ \m{\epsilon(k)=k}.

\noindent - Pour tout entier $k$ et tout \ \m{x\in\D}, on a
\ \m{\epsilon(x+k)=\epsilon(x)+k}.

\noindent - Pour tous entiers $p$, $q$, avec \ \m{q\geq 0}, on a
$$\epsilon(\q{2p+1}{2^{q+1}}) \ = \ 
\epsilon(\q{p}{2^q})\times\epsilon(\q{p+1}{2^q}),$$
o\`u $\times$ est la loi de composition suivante :
$$\alpha\times\beta\ = \ \q{\alpha+\beta}{2}+\q{\Delta_\alpha-\Delta_\beta}
{3+\alpha-\beta}.$$
Cette relation signifie simplement que
$$\chi(E_{\alpha\times\beta},E_\alpha) \ = \ 
\chi(E_\beta,E_{\alpha\times\beta}) \ = \ 0.$$

\bigskip

La construction des pentes des fibr\'es exceptionnels comprises entre $-1$ et
$0$ se fait donc en partant des pentes $-1$ et $0$, correspondant aux
fibr\'es exceptionnels \m{\O(-1)} et $\O$.

On appelle {\em triades} les triplets de fibr\'es exceptionnels de la forme

\noindent\m{(\O(k),\O(k+1),\O(k+2))}, \m{(E_\alpha,E_{\alpha\times\beta},E_\beta)},
\m{(E_{\alpha\times\beta},E_{\beta},E_{\alpha+3})} ou 
\m{(E_{\beta-3},E_\alpha,E_{\alpha\times\beta})}, \m{\alpha} et \m{\beta}
\'etant des \'el\'ements de $\P$ de la forme 
$$\alpha \ = \ \epsilon(\q{p}{2^q}), \ \ \ \
\beta \ = \ \epsilon(\q{p+1}{2^q}).$$
o\`u $p$ et $q$ sont deux entiers avec \ \m{q\geq 0}. Les triades sont
exactement les {\em bases d'h\'elice} de \cite{go_ru}.

On donne maintenant la construction des triades de fibr\'es exceptionnels dont
les pentes sont comprises entre \m{-1} et $0$. Ces triades sont du type
\m{(E_\alpha,E_{\alpha\times\beta},E_\beta)}. La construction
 se fait de la fa\c con suivante, par
r\'ecurrence : on part de la triade \m{(\O(-1),Q^*,\O)}, o\`u $Q$ est le
fibr\'e exceptionnel quotient du morphisme canonique \ 
\m{\O(-1)\lra\O\ot H^0(\O(1))^*}. Supposons la triade \m{(E,F,G)} construite.
Alors on construit les {\it triades adjacentes} \m{(E,H,F)} et \m{(F,K,G)}.
Le fibr\'e $H$ est le noyau du morphisme canonique surjectif
$$F\ot\Hom(F,G)\lra G$$
et $K$ est le conoyau du morphisme canonique injectif
$$E\lra F\ot\Hom(E,F)^*.$$
De plus, le morphisme canonique 
$$E\ot\Hom(E,H)\lra H \ \ {\rm \ \ \ (resp. \ } K\lra G\ot\Hom(K,G)^*
{\rm \ )}$$ 
est surjectif (resp. injectif) et son noyau (resp. conoyau) est isomorphe \`a
\m{G(-3)} (resp. \m{E(3)}). 

\subsection{Suite spectrale de Beilinson g\'en\'eralis\'ee}

A toute triade \m{(E,G,F)} et \`a tout faisceau coh\'erent $\E$ sur \proj{2}
on associe une suite spectrale \m{E^{p,q}_r} de faisceaux coh\'erents sur
\proj{2}, convergeant vers $\E$ en degr\'e 0 et vers 0 en tout autre degr\'e.
Les termes \m{E^{p,q}_1} \'eventuellement non nuls sont
$$E^{-2,q}_1\simeq H^q(\E\ot E^*(-3))\ot E, \ \ 
E^{-1,q}_1\simeq H^q(\E\ot S^*)\ot G, \ \ 
E^{0,q}_1\simeq H^q(\E\ot F^*)\ot F,$$ 
$S$ d\'esignant le fibr\'e exceptionnel conoyau du morphisme canonique injectif

\noindent\m{G\lra F\ot\Hom(G,F)}.

\subsection{S\'erie exceptionnelle associ\'ee \`a un fibr\'e exceptionnel}

Soit $F$ un fibr\'e exceptionnel. Les triades comportant $F$ comme terme
de droite sont de la forme \m{(G_n,G_{n+1},F)}, o\`u la suite de fibr\'es
exceptionnels \m{(G_n)} est enti\`erement d\'etermin\'ee par deux de ses
termes cons\'ecutifs, par exemple \m{G_0} et \m{G_1}, par les suites exactes
$$0\lra G_{n-1}\lra (G_n\ot\Hom(G_{n-1},G_n)^*)\simeq (G_n\ot\Hom(G_n,G_{n+1}))
\lra G_{n+1}\lra 0.$$
On appelle \m{(G_n)} la {\it s\'erie exceptionnelle} \`a gauche associ\'ee
\`a $F$. 
Les couples \m{(\mu(G_n),\Delta(G_n))} sont situ\'es sur la conique 
d'\'equation
$$\Delta \ = \ P(\mu(F)-\mu)-\Delta(F),$$
(ce qui traduit le fait que \ \m{\chi(F,G_n)=0}).

\bigskip

\bigskip

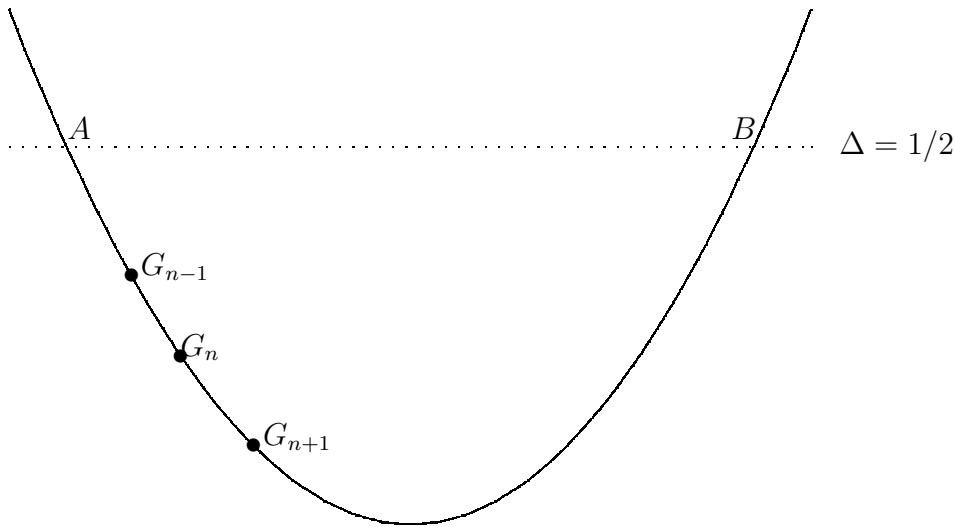
\setlength{\unitlength}{0.240900pt}
\ifx\plotpoint\undefined\newsavebox{\plotpoint}\fi
\begin{picture}(1500,900)(0,0)
\font\gnuplot=cmr10 at 10pt
\gnuplot
\sbox{\plotpoint}{\rule[-0.200pt]{0.400pt}{0.400pt}}%
\put(176,877){\usebox{\plotpoint}}
\multiput(176.58,872.50)(0.493,-1.250){23}{\rule{0.119pt}{1.085pt}}
\multiput(175.17,874.75)(13.000,-29.749){2}{\rule{0.400pt}{0.542pt}}
\multiput(189.58,840.16)(0.492,-1.358){21}{\rule{0.119pt}{1.167pt}}
\multiput(188.17,842.58)(12.000,-29.579){2}{\rule{0.400pt}{0.583pt}}
\multiput(201.58,808.63)(0.493,-1.210){23}{\rule{0.119pt}{1.054pt}}
\multiput(200.17,810.81)(13.000,-28.813){2}{\rule{0.400pt}{0.527pt}}
\multiput(214.58,777.75)(0.493,-1.171){23}{\rule{0.119pt}{1.023pt}}
\multiput(213.17,779.88)(13.000,-27.877){2}{\rule{0.400pt}{0.512pt}}
\multiput(227.58,747.75)(0.493,-1.171){23}{\rule{0.119pt}{1.023pt}}
\multiput(226.17,749.88)(13.000,-27.877){2}{\rule{0.400pt}{0.512pt}}
\multiput(240.58,717.57)(0.492,-1.229){21}{\rule{0.119pt}{1.067pt}}
\multiput(239.17,719.79)(12.000,-26.786){2}{\rule{0.400pt}{0.533pt}}
\multiput(252.58,688.88)(0.493,-1.131){23}{\rule{0.119pt}{0.992pt}}
\multiput(251.17,690.94)(13.000,-26.940){2}{\rule{0.400pt}{0.496pt}}
\multiput(265.58,660.14)(0.493,-1.052){23}{\rule{0.119pt}{0.931pt}}
\put(268,673){$A$}
\multiput(264.17,662.07)(13.000,-25.068){2}{\rule{0.400pt}{0.465pt}}
\multiput(278.58,633.14)(0.493,-1.052){23}{\rule{0.119pt}{0.931pt}}
\multiput(277.17,635.07)(13.000,-25.068){2}{\rule{0.400pt}{0.465pt}}
\multiput(291.58,605.85)(0.492,-1.142){21}{\rule{0.119pt}{1.000pt}}
\multiput(290.17,607.92)(12.000,-24.924){2}{\rule{0.400pt}{0.500pt}}
\multiput(303.58,579.26)(0.493,-1.012){23}{\rule{0.119pt}{0.900pt}}
\multiput(302.17,581.13)(13.000,-24.132){2}{\rule{0.400pt}{0.450pt}}
\multiput(316.58,553.39)(0.493,-0.972){23}{\rule{0.119pt}{0.869pt}}
\multiput(315.17,555.20)(13.000,-23.196){2}{\rule{0.400pt}{0.435pt}}
\multiput(329.58,528.26)(0.492,-1.013){21}{\rule{0.119pt}{0.900pt}}
\multiput(328.17,530.13)(12.000,-22.132){2}{\rule{0.400pt}{0.450pt}}
\multiput(341.58,504.52)(0.493,-0.933){23}{\rule{0.119pt}{0.838pt}}
\multiput(340.17,506.26)(13.000,-22.260){2}{\rule{0.400pt}{0.419pt}}
\multiput(354.58,480.65)(0.493,-0.893){23}{\rule{0.119pt}{0.808pt}}
\multiput(353.17,482.32)(13.000,-21.324){2}{\rule{0.400pt}{0.404pt}}
\multiput(367.58,457.65)(0.493,-0.893){23}{\rule{0.119pt}{0.808pt}}
\multiput(366.17,459.32)(13.000,-21.324){2}{\rule{0.400pt}{0.404pt}}
\put(368,460){\circle*{20}}
\put(383,460){$G_{n-1}$}
\multiput(380.58,434.68)(0.492,-0.884){21}{\rule{0.119pt}{0.800pt}}
\multiput(379.17,436.34)(12.000,-19.340){2}{\rule{0.400pt}{0.400pt}}
\multiput(392.58,413.90)(0.493,-0.814){23}{\rule{0.119pt}{0.746pt}}
\multiput(391.17,415.45)(13.000,-19.451){2}{\rule{0.400pt}{0.373pt}}
\multiput(405.58,392.90)(0.493,-0.814){23}{\rule{0.119pt}{0.746pt}}
\multiput(404.17,394.45)(13.000,-19.451){2}{\rule{0.400pt}{0.373pt}}
\multiput(418.58,372.03)(0.493,-0.774){23}{\rule{0.119pt}{0.715pt}}
\multiput(417.17,373.52)(13.000,-18.515){2}{\rule{0.400pt}{0.358pt}}
\multiput(431.58,351.96)(0.492,-0.798){21}{\rule{0.119pt}{0.733pt}}
\multiput(430.17,353.48)(12.000,-17.478){2}{\rule{0.400pt}{0.367pt}}
\multiput(443.58,333.29)(0.493,-0.695){23}{\rule{0.119pt}{0.654pt}}
\put(445,333){\circle*{20}}
\put(445,333){$G_{n}$}
\multiput(442.17,334.64)(13.000,-16.643){2}{\rule{0.400pt}{0.327pt}}
\multiput(456.58,315.29)(0.493,-0.695){23}{\rule{0.119pt}{0.654pt}}
\multiput(455.17,316.64)(13.000,-16.643){2}{\rule{0.400pt}{0.327pt}}
\multiput(469.58,297.23)(0.492,-0.712){21}{\rule{0.119pt}{0.667pt}}
\multiput(468.17,298.62)(12.000,-15.616){2}{\rule{0.400pt}{0.333pt}}
\multiput(481.58,280.41)(0.493,-0.655){23}{\rule{0.119pt}{0.623pt}}
\multiput(480.17,281.71)(13.000,-15.707){2}{\rule{0.400pt}{0.312pt}}
\multiput(494.58,263.54)(0.493,-0.616){23}{\rule{0.119pt}{0.592pt}}
\multiput(493.17,264.77)(13.000,-14.771){2}{\rule{0.400pt}{0.296pt}}
\multiput(507.58,247.67)(0.493,-0.576){23}{\rule{0.119pt}{0.562pt}}
\multiput(506.17,248.83)(13.000,-13.834){2}{\rule{0.400pt}{0.281pt}}
\multiput(520.58,232.65)(0.492,-0.582){21}{\rule{0.119pt}{0.567pt}}
\multiput(519.17,233.82)(12.000,-12.824){2}{\rule{0.400pt}{0.283pt}}
\multiput(532.58,218.80)(0.493,-0.536){23}{\rule{0.119pt}{0.531pt}}
\multiput(531.17,219.90)(13.000,-12.898){2}{\rule{0.400pt}{0.265pt}}
\multiput(545.58,204.80)(0.493,-0.536){23}{\rule{0.119pt}{0.531pt}}
\multiput(544.17,205.90)(13.000,-12.898){2}{\rule{0.400pt}{0.265pt}}
\multiput(558.00,191.92)(0.539,-0.492){21}{\rule{0.533pt}{0.119pt}}
\multiput(558.00,192.17)(11.893,-12.000){2}{\rule{0.267pt}{0.400pt}}
\put(560,192){\circle*{20}}
\put(575,192){$G_{n+1}$}
\multiput(571.00,179.92)(0.496,-0.492){21}{\rule{0.500pt}{0.119pt}}
\multiput(571.00,180.17)(10.962,-12.000){2}{\rule{0.250pt}{0.400pt}}
\multiput(583.00,167.92)(0.590,-0.492){19}{\rule{0.573pt}{0.118pt}}
\multiput(583.00,168.17)(11.811,-11.000){2}{\rule{0.286pt}{0.400pt}}
\multiput(596.00,156.92)(0.590,-0.492){19}{\rule{0.573pt}{0.118pt}}
\multiput(596.00,157.17)(11.811,-11.000){2}{\rule{0.286pt}{0.400pt}}
\multiput(609.00,145.92)(0.600,-0.491){17}{\rule{0.580pt}{0.118pt}}
\multiput(609.00,146.17)(10.796,-10.000){2}{\rule{0.290pt}{0.400pt}}
\multiput(621.00,135.93)(0.728,-0.489){15}{\rule{0.678pt}{0.118pt}}
\multiput(621.00,136.17)(11.593,-9.000){2}{\rule{0.339pt}{0.400pt}}
\multiput(634.00,126.93)(0.824,-0.488){13}{\rule{0.750pt}{0.117pt}}
\multiput(634.00,127.17)(11.443,-8.000){2}{\rule{0.375pt}{0.400pt}}
\multiput(647.00,118.93)(0.824,-0.488){13}{\rule{0.750pt}{0.117pt}}
\multiput(647.00,119.17)(11.443,-8.000){2}{\rule{0.375pt}{0.400pt}}
\multiput(660.00,110.93)(0.758,-0.488){13}{\rule{0.700pt}{0.117pt}}
\multiput(660.00,111.17)(10.547,-8.000){2}{\rule{0.350pt}{0.400pt}}
\multiput(672.00,102.93)(1.123,-0.482){9}{\rule{0.967pt}{0.116pt}}
\multiput(672.00,103.17)(10.994,-6.000){2}{\rule{0.483pt}{0.400pt}}
\multiput(685.00,96.93)(1.123,-0.482){9}{\rule{0.967pt}{0.116pt}}
\multiput(685.00,97.17)(10.994,-6.000){2}{\rule{0.483pt}{0.400pt}}
\multiput(698.00,90.93)(1.123,-0.482){9}{\rule{0.967pt}{0.116pt}}
\multiput(698.00,91.17)(10.994,-6.000){2}{\rule{0.483pt}{0.400pt}}
\multiput(711.00,84.94)(1.651,-0.468){5}{\rule{1.300pt}{0.113pt}}
\multiput(711.00,85.17)(9.302,-4.000){2}{\rule{0.650pt}{0.400pt}}
\multiput(723.00,80.94)(1.797,-0.468){5}{\rule{1.400pt}{0.113pt}}
\multiput(723.00,81.17)(10.094,-4.000){2}{\rule{0.700pt}{0.400pt}}
\multiput(736.00,76.95)(2.695,-0.447){3}{\rule{1.833pt}{0.108pt}}
\multiput(736.00,77.17)(9.195,-3.000){2}{\rule{0.917pt}{0.400pt}}
\multiput(749.00,73.95)(2.472,-0.447){3}{\rule{1.700pt}{0.108pt}}
\multiput(749.00,74.17)(8.472,-3.000){2}{\rule{0.850pt}{0.400pt}}
\put(761,70.17){\rule{2.700pt}{0.400pt}}
\multiput(761.00,71.17)(7.396,-2.000){2}{\rule{1.350pt}{0.400pt}}
\put(774,68.67){\rule{3.132pt}{0.400pt}}
\multiput(774.00,69.17)(6.500,-1.000){2}{\rule{1.566pt}{0.400pt}}
\put(787,67.67){\rule{3.132pt}{0.400pt}}
\multiput(787.00,68.17)(6.500,-1.000){2}{\rule{1.566pt}{0.400pt}}
\put(812,67.67){\rule{3.132pt}{0.400pt}}
\multiput(812.00,67.17)(6.500,1.000){2}{\rule{1.566pt}{0.400pt}}
\put(825,68.67){\rule{3.132pt}{0.400pt}}
\multiput(825.00,68.17)(6.500,1.000){2}{\rule{1.566pt}{0.400pt}}
\put(838,70.17){\rule{2.700pt}{0.400pt}}
\multiput(838.00,69.17)(7.396,2.000){2}{\rule{1.350pt}{0.400pt}}
\multiput(851.00,72.61)(2.472,0.447){3}{\rule{1.700pt}{0.108pt}}
\multiput(851.00,71.17)(8.472,3.000){2}{\rule{0.850pt}{0.400pt}}
\multiput(863.00,75.61)(2.695,0.447){3}{\rule{1.833pt}{0.108pt}}
\multiput(863.00,74.17)(9.195,3.000){2}{\rule{0.917pt}{0.400pt}}
\multiput(876.00,78.60)(1.797,0.468){5}{\rule{1.400pt}{0.113pt}}
\multiput(876.00,77.17)(10.094,4.000){2}{\rule{0.700pt}{0.400pt}}
\multiput(889.00,82.60)(1.651,0.468){5}{\rule{1.300pt}{0.113pt}}
\multiput(889.00,81.17)(9.302,4.000){2}{\rule{0.650pt}{0.400pt}}
\multiput(901.00,86.59)(1.123,0.482){9}{\rule{0.967pt}{0.116pt}}
\multiput(901.00,85.17)(10.994,6.000){2}{\rule{0.483pt}{0.400pt}}
\multiput(914.00,92.59)(1.123,0.482){9}{\rule{0.967pt}{0.116pt}}
\multiput(914.00,91.17)(10.994,6.000){2}{\rule{0.483pt}{0.400pt}}
\multiput(927.00,98.59)(1.123,0.482){9}{\rule{0.967pt}{0.116pt}}
\multiput(927.00,97.17)(10.994,6.000){2}{\rule{0.483pt}{0.400pt}}
\multiput(940.00,104.59)(0.758,0.488){13}{\rule{0.700pt}{0.117pt}}
\multiput(940.00,103.17)(10.547,8.000){2}{\rule{0.350pt}{0.400pt}}
\multiput(952.00,112.59)(0.824,0.488){13}{\rule{0.750pt}{0.117pt}}
\multiput(952.00,111.17)(11.443,8.000){2}{\rule{0.375pt}{0.400pt}}
\multiput(965.00,120.59)(0.824,0.488){13}{\rule{0.750pt}{0.117pt}}
\multiput(965.00,119.17)(11.443,8.000){2}{\rule{0.375pt}{0.400pt}}
\multiput(978.00,128.59)(0.728,0.489){15}{\rule{0.678pt}{0.118pt}}
\multiput(978.00,127.17)(11.593,9.000){2}{\rule{0.339pt}{0.400pt}}
\multiput(991.00,137.58)(0.600,0.491){17}{\rule{0.580pt}{0.118pt}}
\multiput(991.00,136.17)(10.796,10.000){2}{\rule{0.290pt}{0.400pt}}
\multiput(1003.00,147.58)(0.590,0.492){19}{\rule{0.573pt}{0.118pt}}
\multiput(1003.00,146.17)(11.811,11.000){2}{\rule{0.286pt}{0.400pt}}
\multiput(1016.00,158.58)(0.590,0.492){19}{\rule{0.573pt}{0.118pt}}
\multiput(1016.00,157.17)(11.811,11.000){2}{\rule{0.286pt}{0.400pt}}
\multiput(1029.00,169.58)(0.496,0.492){21}{\rule{0.500pt}{0.119pt}}
\multiput(1029.00,168.17)(10.962,12.000){2}{\rule{0.250pt}{0.400pt}}
\multiput(1041.00,181.58)(0.539,0.492){21}{\rule{0.533pt}{0.119pt}}
\multiput(1041.00,180.17)(11.893,12.000){2}{\rule{0.267pt}{0.400pt}}
\multiput(1054.58,193.00)(0.493,0.536){23}{\rule{0.119pt}{0.531pt}}
\multiput(1053.17,193.00)(13.000,12.898){2}{\rule{0.400pt}{0.265pt}}
\multiput(1067.58,207.00)(0.493,0.536){23}{\rule{0.119pt}{0.531pt}}
\multiput(1066.17,207.00)(13.000,12.898){2}{\rule{0.400pt}{0.265pt}}
\multiput(1080.58,221.00)(0.492,0.582){21}{\rule{0.119pt}{0.567pt}}
\multiput(1079.17,221.00)(12.000,12.824){2}{\rule{0.400pt}{0.283pt}}
\multiput(1092.58,235.00)(0.493,0.576){23}{\rule{0.119pt}{0.562pt}}
\multiput(1091.17,235.00)(13.000,13.834){2}{\rule{0.400pt}{0.281pt}}
\multiput(1105.58,250.00)(0.493,0.616){23}{\rule{0.119pt}{0.592pt}}
\multiput(1104.17,250.00)(13.000,14.771){2}{\rule{0.400pt}{0.296pt}}
\multiput(1118.58,266.00)(0.493,0.655){23}{\rule{0.119pt}{0.623pt}}
\multiput(1117.17,266.00)(13.000,15.707){2}{\rule{0.400pt}{0.312pt}}
\multiput(1131.58,283.00)(0.492,0.712){21}{\rule{0.119pt}{0.667pt}}
\multiput(1130.17,283.00)(12.000,15.616){2}{\rule{0.400pt}{0.333pt}}
\multiput(1143.58,300.00)(0.493,0.695){23}{\rule{0.119pt}{0.654pt}}
\multiput(1142.17,300.00)(13.000,16.643){2}{\rule{0.400pt}{0.327pt}}
\multiput(1156.58,318.00)(0.493,0.695){23}{\rule{0.119pt}{0.654pt}}
\multiput(1155.17,318.00)(13.000,16.643){2}{\rule{0.400pt}{0.327pt}}
\multiput(1169.58,336.00)(0.492,0.798){21}{\rule{0.119pt}{0.733pt}}
\multiput(1168.17,336.00)(12.000,17.478){2}{\rule{0.400pt}{0.367pt}}
\multiput(1181.58,355.00)(0.493,0.774){23}{\rule{0.119pt}{0.715pt}}
\multiput(1180.17,355.00)(13.000,18.515){2}{\rule{0.400pt}{0.358pt}}
\multiput(1194.58,375.00)(0.493,0.814){23}{\rule{0.119pt}{0.746pt}}
\multiput(1193.17,375.00)(13.000,19.451){2}{\rule{0.400pt}{0.373pt}}
\multiput(1207.58,396.00)(0.493,0.814){23}{\rule{0.119pt}{0.746pt}}
\multiput(1206.17,396.00)(13.000,19.451){2}{\rule{0.400pt}{0.373pt}}
\multiput(1220.58,417.00)(0.492,0.884){21}{\rule{0.119pt}{0.800pt}}
\multiput(1219.17,417.00)(12.000,19.340){2}{\rule{0.400pt}{0.400pt}}
\multiput(1232.58,438.00)(0.493,0.893){23}{\rule{0.119pt}{0.808pt}}
\multiput(1231.17,438.00)(13.000,21.324){2}{\rule{0.400pt}{0.404pt}}
\multiput(1245.58,461.00)(0.493,0.893){23}{\rule{0.119pt}{0.808pt}}
\multiput(1244.17,461.00)(13.000,21.324){2}{\rule{0.400pt}{0.404pt}}
\multiput(1258.58,484.00)(0.493,0.933){23}{\rule{0.119pt}{0.838pt}}
\multiput(1257.17,484.00)(13.000,22.260){2}{\rule{0.400pt}{0.419pt}}
\multiput(1271.58,508.00)(0.492,1.013){21}{\rule{0.119pt}{0.900pt}}
\multiput(1270.17,508.00)(12.000,22.132){2}{\rule{0.400pt}{0.450pt}}
\multiput(1283.58,532.00)(0.493,0.972){23}{\rule{0.119pt}{0.869pt}}
\multiput(1282.17,532.00)(13.000,23.196){2}{\rule{0.400pt}{0.435pt}}
\multiput(1296.58,557.00)(0.493,1.012){23}{\rule{0.119pt}{0.900pt}}
\multiput(1295.17,557.00)(13.000,24.132){2}{\rule{0.400pt}{0.450pt}}
\multiput(1309.58,583.00)(0.492,1.142){21}{\rule{0.119pt}{1.000pt}}
\multiput(1308.17,583.00)(12.000,24.924){2}{\rule{0.400pt}{0.500pt}}
\multiput(1321.58,610.00)(0.493,1.052){23}{\rule{0.119pt}{0.931pt}}
\multiput(1320.17,610.00)(13.000,25.068){2}{\rule{0.400pt}{0.465pt}}
\multiput(1334.58,637.00)(0.493,1.052){23}{\rule{0.119pt}{0.931pt}}
\multiput(1333.17,637.00)(13.000,25.068){2}{\rule{0.400pt}{0.465pt}}
\multiput(1347.58,664.00)(0.493,1.131){23}{\rule{0.119pt}{0.992pt}}
\multiput(1346.17,664.00)(13.000,26.940){2}{\rule{0.400pt}{0.496pt}}
\put(1310,673){$B$}
\multiput(1360.58,693.00)(0.492,1.229){21}{\rule{0.119pt}{1.067pt}}
\multiput(1359.17,693.00)(12.000,26.786){2}{\rule{0.400pt}{0.533pt}}
\multiput(1372.58,722.00)(0.493,1.171){23}{\rule{0.119pt}{1.023pt}}
\multiput(1371.17,722.00)(13.000,27.877){2}{\rule{0.400pt}{0.512pt}}
\multiput(1385.58,752.00)(0.493,1.171){23}{\rule{0.119pt}{1.023pt}}
\multiput(1384.17,752.00)(13.000,27.877){2}{\rule{0.400pt}{0.512pt}}
\multiput(1398.58,782.00)(0.493,1.210){23}{\rule{0.119pt}{1.054pt}}
\multiput(1397.17,782.00)(13.000,28.813){2}{\rule{0.400pt}{0.527pt}}
\multiput(1411.58,813.00)(0.492,1.358){21}{\rule{0.119pt}{1.167pt}}
\multiput(1410.17,813.00)(12.000,29.579){2}{\rule{0.400pt}{0.583pt}}
\multiput(1423.58,845.00)(0.493,1.250){23}{\rule{0.119pt}{1.085pt}}
\multiput(1422.17,845.00)(13.000,29.749){2}{\rule{0.400pt}{0.542pt}}
\put(800.0,68.0){\rule[-0.200pt]{2.891pt}{0.400pt}}
\put(176,661){\usebox{\plotpoint}}
\put(176.00,661.00){\usebox{\plotpoint}}
\put(196.76,661.00){\usebox{\plotpoint}}
\multiput(201,661)(20.756,0.000){0}{\usebox{\plotpoint}}
\put(217.51,661.00){\usebox{\plotpoint}}
\put(238.27,661.00){\usebox{\plotpoint}}
\multiput(240,661)(20.756,0.000){0}{\usebox{\plotpoint}}
\put(259.02,661.00){\usebox{\plotpoint}}
\multiput(265,661)(20.756,0.000){0}{\usebox{\plotpoint}}
\put(279.78,661.00){\usebox{\plotpoint}}
\put(300.53,661.00){\usebox{\plotpoint}}
\multiput(303,661)(20.756,0.000){0}{\usebox{\plotpoint}}
\put(321.29,661.00){\usebox{\plotpoint}}
\multiput(329,661)(20.756,0.000){0}{\usebox{\plotpoint}}
\put(342.04,661.00){\usebox{\plotpoint}}
\put(362.80,661.00){\usebox{\plotpoint}}
\multiput(367,661)(20.756,0.000){0}{\usebox{\plotpoint}}
\put(383.55,661.00){\usebox{\plotpoint}}
\put(404.31,661.00){\usebox{\plotpoint}}
\multiput(405,661)(20.756,0.000){0}{\usebox{\plotpoint}}
\put(425.07,661.00){\usebox{\plotpoint}}
\multiput(431,661)(20.756,0.000){0}{\usebox{\plotpoint}}
\put(445.82,661.00){\usebox{\plotpoint}}
\put(466.58,661.00){\usebox{\plotpoint}}
\multiput(469,661)(20.756,0.000){0}{\usebox{\plotpoint}}
\put(487.33,661.00){\usebox{\plotpoint}}
\multiput(494,661)(20.756,0.000){0}{\usebox{\plotpoint}}
\put(508.09,661.00){\usebox{\plotpoint}}
\put(528.84,661.00){\usebox{\plotpoint}}
\multiput(532,661)(20.756,0.000){0}{\usebox{\plotpoint}}
\put(549.60,661.00){\usebox{\plotpoint}}
\put(570.35,661.00){\usebox{\plotpoint}}
\multiput(571,661)(20.756,0.000){0}{\usebox{\plotpoint}}
\put(591.11,661.00){\usebox{\plotpoint}}
\multiput(596,661)(20.756,0.000){0}{\usebox{\plotpoint}}
\put(611.87,661.00){\usebox{\plotpoint}}
\put(632.62,661.00){\usebox{\plotpoint}}
\multiput(634,661)(20.756,0.000){0}{\usebox{\plotpoint}}
\put(653.38,661.00){\usebox{\plotpoint}}
\multiput(660,661)(20.756,0.000){0}{\usebox{\plotpoint}}
\put(674.13,661.00){\usebox{\plotpoint}}
\put(694.89,661.00){\usebox{\plotpoint}}
\multiput(698,661)(20.756,0.000){0}{\usebox{\plotpoint}}
\put(715.64,661.00){\usebox{\plotpoint}}
\multiput(723,661)(20.756,0.000){0}{\usebox{\plotpoint}}
\put(736.40,661.00){\usebox{\plotpoint}}
\put(757.15,661.00){\usebox{\plotpoint}}
\multiput(761,661)(20.756,0.000){0}{\usebox{\plotpoint}}
\put(777.91,661.00){\usebox{\plotpoint}}
\put(798.66,661.00){\usebox{\plotpoint}}
\multiput(800,661)(20.756,0.000){0}{\usebox{\plotpoint}}
\put(819.42,661.00){\usebox{\plotpoint}}
\multiput(825,661)(20.756,0.000){0}{\usebox{\plotpoint}}
\put(840.18,661.00){\usebox{\plotpoint}}
\put(860.93,661.00){\usebox{\plotpoint}}
\multiput(863,661)(20.756,0.000){0}{\usebox{\plotpoint}}
\put(881.69,661.00){\usebox{\plotpoint}}
\multiput(889,661)(20.756,0.000){0}{\usebox{\plotpoint}}
\put(902.44,661.00){\usebox{\plotpoint}}
\put(923.20,661.00){\usebox{\plotpoint}}
\multiput(927,661)(20.756,0.000){0}{\usebox{\plotpoint}}
\put(943.95,661.00){\usebox{\plotpoint}}
\put(964.71,661.00){\usebox{\plotpoint}}
\multiput(965,661)(20.756,0.000){0}{\usebox{\plotpoint}}
\put(985.46,661.00){\usebox{\plotpoint}}
\multiput(991,661)(20.756,0.000){0}{\usebox{\plotpoint}}
\put(1006.22,661.00){\usebox{\plotpoint}}
\put(1026.98,661.00){\usebox{\plotpoint}}
\multiput(1029,661)(20.756,0.000){0}{\usebox{\plotpoint}}
\put(1047.73,661.00){\usebox{\plotpoint}}
\multiput(1054,661)(20.756,0.000){0}{\usebox{\plotpoint}}
\put(1068.49,661.00){\usebox{\plotpoint}}
\put(1089.24,661.00){\usebox{\plotpoint}}
\multiput(1092,661)(20.756,0.000){0}{\usebox{\plotpoint}}
\put(1110.00,661.00){\usebox{\plotpoint}}
\put(1130.75,661.00){\usebox{\plotpoint}}
\multiput(1131,661)(20.756,0.000){0}{\usebox{\plotpoint}}
\put(1151.51,661.00){\usebox{\plotpoint}}
\multiput(1156,661)(20.756,0.000){0}{\usebox{\plotpoint}}
\put(1172.26,661.00){\usebox{\plotpoint}}
\put(1193.02,661.00){\usebox{\plotpoint}}
\multiput(1194,661)(20.756,0.000){0}{\usebox{\plotpoint}}
\put(1213.77,661.00){\usebox{\plotpoint}}
\multiput(1220,661)(20.756,0.000){0}{\usebox{\plotpoint}}
\put(1234.53,661.00){\usebox{\plotpoint}}
\put(1255.29,661.00){\usebox{\plotpoint}}
\multiput(1258,661)(20.756,0.000){0}{\usebox{\plotpoint}}
\put(1276.04,661.00){\usebox{\plotpoint}}
\multiput(1283,661)(20.756,0.000){0}{\usebox{\plotpoint}}
\put(1296.80,661.00){\usebox{\plotpoint}}
\put(1317.55,661.00){\usebox{\plotpoint}}
\multiput(1321,661)(20.756,0.000){0}{\usebox{\plotpoint}}
\put(1338.31,661.00){\usebox{\plotpoint}}
\put(1359.06,661.00){\usebox{\plotpoint}}
\multiput(1360,661)(20.756,0.000){0}{\usebox{\plotpoint}}
\put(1379.82,661.00){\usebox{\plotpoint}}
\multiput(1385,661)(20.756,0.000){0}{\usebox{\plotpoint}}
\put(1400.57,661.00){\usebox{\plotpoint}}
\put(1421.33,661.00){\usebox{\plotpoint}}
\multiput(1423,661)(20.756,0.000){0}{\usebox{\plotpoint}}
\put(1436,661){\usebox{\plotpoint}}
\put(1480,650){$\Delta=1/2$}
\end{picture}

\bigskip

Dans la figure ci-dessus, les points $A$ et $B$ sont les intersections de
cette conique avec la droite d'\'equation \ \m{\Delta=1/2}. On a
$$\lim_{n\rightarrow -\infty}=A \ \ \ {\rm et} \ \ \ 
\lim_{n\rightarrow\infty}=B.$$
Remarquons que \ \m{\mu(B)-\mu(A)<3}.

Si \m{F=\O}, il existe une unique paire \m{(G_n,G_{n+1})} telle que \
\m{\mu(G_{n+1})-\mu(G_n)\geq 1}, c'est \m{(\O(-2),\O(-1))}. Supposons que
\ \m{-1<\mu(F)<0}. Il existe alors une unique triade de la forme 
\m{(E,F,G)}, avec \ \m{-1\leq\mu(E)<\mu(G)\leq 0}. On en d\'eduit que
\m{(G(-3),E)} est une des paires \m{(G_n,G_{n+1})}. On peut supposer que
\ \m{(G(-3),E)=(G_0,G_1)}. On a \ \m{\mu(G_1)-\mu(G_0)\geq 2}, et \m{(G_0,G_1)}
est l'unique paire  \m{(G_n,G_{n+1})} telle que 
 
\noindent\m{\mu(G_{n+1})-\mu(G_n)\geq 1}. On l'appelle la paire {\it initiale}
de la s\'erie \m{(G_n)}.

\bigskip

\begin{xlemm}
Le fibr\'e vectoriel \m{G_0^*\ot G_1} est engendr\'e par ses sections globales. 
\end{xlemm}

\dem D'apr\`es la construction de \m{(G_0,G_1)}, il suffit de prouver le
r\'esultat \hbox{suivant :} si \m{(A,B,C)} est une triade de fibr\'es
exceptionnels telle que \ \m{\mu(C)-\mu(A)\leq 1}, les fibr\'es 
\m{B^*\ot A(3)}, \m{C^*\ot B(3)} et \m{C^*\ot A(3)} sont engendr\'es par leurs
sections globales. On d\'emontre cela par r\'ecurrence : il faut montrer que 
si c'est vrai pour une triade, c'est vrai pour les deux triades adjacentes.
Supposons que ce soit vrai pour \m{(A,B,C)}. Soient $H$ le noyau du morphisme
canonique surjectif
$$B\ot\Hom(B,C)\lra C$$
et K le conoyau du morphisme canonique injectif
$$A\lra B\ot\Hom(A,B)^*.$$
Il faut montrer que le r\'esultat est vrai pour les triades \m{(A,H,B)} et
\m{(B,K,C)}. En consid\'erant la triade {\it duale}
\m{(C^*(-1),B^*(-1),A^*(-1))}, on voit qu'il suffit de consid\'erer 
\m{(A,H,B)}. On a une suite exacte
$$0\lra H\lra B\ot\Hom(B,C)\lra C\lra 0.$$
On en d\'eduit un morphisme surjectif
$$B^*(3)\ot A\ot\Hom(B,C)^*\lra H^*(3)\ot A.$$
Puisque \m{B^*(3)\ot A} est engendr\'e par ses sections globales (hypoth\`ese
de r\'ecurrence), il en est de m\^eme de \m{H^*(3)\ot A}. On a d'autre
part une suite exacte
$$0\lra C(-3)\lra A\ot\Hom(C(-3),A)^*\lra H\lra 0,$$
d'o\`u on d\'eduit un morphisme surjectif
$$B^*(3)\ot A\ot\Hom(C(-3),A)^*\lra B^*(3)\ot H,$$
d'o\`u on d\'eduit que \m{B^*(3)\ot H} est engendr\'e par ses sections
globales. \carre

\bigskip

\begin{xlemm}
Pour tout entier $n$, on a \m{n\geq 1} si et seulement si pour tous entiers
$a$, $b$, $c$ positifs ou nuls, le fibr\'e vectoriel
$$(G_n\ot\cx{a})\oplus(G_{n+1}\ot\cx{b})\oplus(F\ot\cx{c})$$
est prioritaire.
\end{xlemm}

\dem Imm\'ediat. \carre

\bigskip

On d\'efinit de m\^eme la {\em s\'erie exceptionnelle \`a droite} \m{(H_n)}
associ\'ee \`a $F$. On a \break \m{H_n=G_n(3)} pour tout $n$.

\subsection{\'Etude de {\bf T}}

L'ensemble {\bf T} est construit comme une union croissante de sous-ensembles
$$T_0=\lbrace(\O(-1),Q^*,\O)\rbrace\subset T_1\subset\ldots T_n\subset
T_{n+1}\subset\ldots$$
$$T=\bigcup_{n\geq 0}T_n,$$
o\`u $T_n$ est l'ensemble des triades \m{(E_\alpha,E_{\alpha\times\beta},
E_\beta)}, $\alpha$, $\beta$ \'etant de la forme
$$\alpha=\epsilon(\q{p}{2^n}), \ \ \beta=\epsilon(\q{p+1}{2^n}),$$
avec $p$ entier. Si $n>0$, les triades de \m{T_n\backslash T_{n-1}} forment
une suite \m{t_0^{(n)}}, \ldots, \m{t_{2^n-1}^{(n)}},
$$t_i^{(n)}\ = \ (E_{\alpha(\q{i}{2^n})},E_{\alpha(\q{2i+1}{2^{n+1}})},
E_{\alpha(\q{i+1}{2^n})}).$$
On a
$$\mu(E_{\alpha(\q{i}{2^n})}) \ < \ \mu(E_{\alpha(\q{2i+1}{2^{n+1}})}) \ < \
\mu(E_{\alpha(\q{i+1}{2^n})}),$$
et dans le plan de coordonn\'ees \m{(\mu,\Delta)},
\m{E_{\alpha(\q{2i+1}{2^{n+1}})}} est situ\'e au dessus de la droite\break
\m{E_{\alpha(\q{i}{2^n})}E_{\alpha(\q{i+1}{2^n})}}.

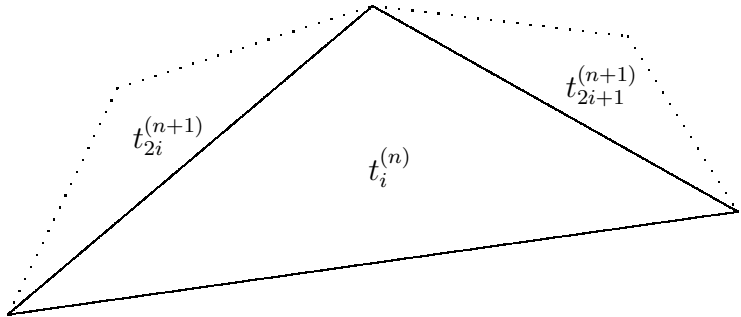
\setlength{\unitlength}{0.240900pt}
\ifx\plotpoint\undefined\newsavebox{\plotpoint}\fi
\begin{picture}(1500,900)(0,0)
\font\gnuplot=cmr10 at 10pt
\gnuplot
\sbox{\plotpoint}{\rule[-0.200pt]{0.400pt}{0.400pt}}%
\put(233,68){\usebox{\plotpoint}}
\multiput(233.00,68.58)(3.543,0.500){321}{\rule{2.930pt}{0.120pt}}
\multiput(233.00,67.17)(1139.919,162.000){2}{\rule{1.465pt}{0.400pt}}
\put(233,68){\usebox{\plotpoint}}
\multiput(233.00,68.58)(0.591,0.500){967}{\rule{0.573pt}{0.120pt}}
\multiput(233.00,67.17)(571.812,485.000){2}{\rule{0.286pt}{0.400pt}}
\put(806,553){\usebox{\plotpoint}}
\multiput(806.00,551.92)(0.887,-0.500){643}{\rule{0.810pt}{0.120pt}}
\multiput(806.00,552.17)(571.320,-323.000){2}{\rule{0.405pt}{0.400pt}}
\put(233,68){\usebox{\plotpoint}}
\multiput(233,68)(9.029,18.689){20}{\usebox{\plotpoint}}
\put(405,424){\usebox{\plotpoint}}
\put(405,424){\usebox{\plotpoint}}
\multiput(405,424)(19.758,6.356){21}{\usebox{\plotpoint}}
\put(806,553){\usebox{\plotpoint}}
\put(806,553){\usebox{\plotpoint}}
\multiput(806,553)(20.608,-2.467){20}{\usebox{\plotpoint}}
\put(1207,505){\usebox{\plotpoint}}
\put(1207,505){\usebox{\plotpoint}}
\multiput(1207,505)(11.006,-17.597){16}{\usebox{\plotpoint}}
\put(1379,230){\usebox{\plotpoint}}
\put(800,280){$t_i^{(n)}$}
\put(430,330){$t_{2i}^{(n+1)}$}
\put(1110,410){$t_{2i+1}^{(n+1)}$}
\end{picture}

Le segment de conique \m{E_{-1}E_0} de \m{\T_{(E_{-1},E_{\q{1}{2}},E_0)}}
n'est autre que la courbe \ \m{\Delta=-\q{\mu(\mu+1)}{2}}. On en d\'eduit
imm\'ediatement le

\bigskip

\begin{xlemm}
Soit \ \m{Z = \bigcup_{(E,F,G)\in{\bf T}}\T_{(E,F,G)}}. Alors, si
\ \m{(\mu,\Delta)\in Z}, on a \ \m{(\mu,\Delta')\in Z} \ si
$$-\q{\mu(\mu+1)}{2} \ \leq \ \Delta' \ \leq \ \Delta.$$
\end{xlemm}

\section{Fibr\'es prioritaires g\'en\'eriques}

\subsection{Cohomologie naturelle}

\begin{xlemm}
Soient $F$ un fibr\'e exceptionnel, $r$, \m{c_1}, \m{c_2} des entiers tels que
\m{r\geq 2},\break \m{\mu(F)-x_F<\mu\leq\mu(F)} \ et \ \m{\Delta=\delta(\mu)}.
Alors il existe un fibr\'e vectoriel stable $\E$ de rang $r$ et de classes de
Chern $c_1$, $c_2$, tel que \ \m{\Ext^1(\E,F)=\lbrace 0\rbrace}.
\end{xlemm}

\dem On consid\`ere la suite \m{(G_n)} de fibr\'es exceptionnels du \para 2. 
Soient $n$ un entier et $\E$ un faisceau semi-stable de rang $r$ et de
classes de Chern \m{c_1}, \m{c_2}. On pose
$$k = \chi(\E,F), \ \ \ m_n \ = \ -\chi(\E\ot G_n^*(-3)),$$
qui sont ind\'ependants de $\E$. Ces entiers sont positifs : pour le premier,
cela d\'ecoule du fait que le point correspondant \`a $\E$ est situ\'e sous
la conique donnant l'\'equation de \m{\delta(\mu)} sur 
\m{\rbrack\mu(F),\mu(F)+x_F\lbrack}. Pour le second on utilise le fait que
\m{H^0(\E\ot G_n^*(-3))} et \m{H^2(\E\ot G_n^*(-3))} sont nuls.
On consid\`ere les triades \m{(F,G_{p-1}(3),G_p(3))}. Ceci sugg\`ere de trouver
$\E$ comme noyau d'un morphisme surjectif ad\'equat
$$\theta : 
(F\ot\cx{k})\oplus(G_{p-1}(3)\ot\cx{m_{p+1}})\lra G_p(3)\ot\cx{m_p}.$$
Un tel fibr\'e a en effet les bons rang et classes de Chern, et de plus on
a \ \m{\Ext^1(\E,F)=\lbrace 0\rbrace}. Pour montrer que $\E$ se d\'eforme en
fibr\'e stable, il suffit qu'il soit prioritaire, car le champ des
faisceaux prioritaires est irr\'eductible (cf. \cite{hi_la}). On prend
\ \m{p=1}, c'est-\`a-dire qu'on consid\`ere des morphismes
$$(F\ot\cx{k})\oplus(G_0(3)\ot\cx{m_2})\lra G_1(3)\ot\cx{m_1}.$$
Alors on a \ \m{\mu(G_1(3))-\mu(G_0(3))\geq 1}, donc \m{\mu(G_1(3))-\mu(F) > 1},
et la paire \m{(F,G_1(3))} est initiale dans la s\'erie qui la contient. Ceci
entraine que le faisceau des morphismes pr\'ec\'edents est engendr\'e par
ses sections globales. Comme \ \m{r\geq 2}, il existe un morphisme
$$\theta : (F\ot\cx{k})\oplus(G_0(3)\ot\cx{m_2})\lra G_1(3)\ot\cx{m_1}$$
qui est surjectif. Soit 
$$\E \ = \ \ker(\theta).$$
Il reste \`a montrer que $\E$ est prioritaire, c'est-\`a-dire que
\ \m{\Hom(\E,\E(-2))=\lbrace 0\rbrace}. On a une suite exacte
$$0\lra\E\lra (F\ot\cx{k})\oplus(G_0(3)\ot\cx{m_2})\lra G_1(3)\ot\cx{m_1}
\lra 0,$$
d'o\`u on d\'eduit que
$$\Hom(\E,\E(-2))\ \subset \ (\Hom(\E,F(-2))\ot\cx{k})\oplus
(\Hom(\E,G_0(1))\ot\cx{m_2}).$$
Il faut montrer que
$$\Hom((\E,F(-2))=\Hom(\E,G_0(1))=\lbrace 0\rbrace.$$

Montrons d'abord que \ \m{\Hom((\E,F(-2))=\lbrace 0\rbrace}.
D'apr\`es la suite exacte pr\'ec\'edente, on a une suite exacte
$$(\Hom(F,F(-2))\ot\cx{k})\oplus(\Hom(G_0(3),F(-2))\ot\cx{m_2})\lra
\Hom((\E,F(-2)) \ \ \ \ \ \ \ \ \ \ \ \ \ \ \ \ $$
$$ \ \ \ \ \ \ \ \ \ \ \ \ \ \ \ \ \ \lra\Ext^1(G_1(3),F(-2))\ot\cx{m_1}.$$
On a \ \m{\Hom(F,F(-2))=\Hom(G_0(3),F(-2))=\lbrace 0\rbrace}, car
\ \m{\mu(G_0(3))>\mu(F)>\mu(F(-2))}. D'autre part,
$$\Ext^1(G_1(3),F(-2))\ \simeq\ \Ext^1(F(-2),G_1)^*$$
par dualit\'e de Serre. Pour montrer que \ \m{\Ext^1(F(-2),G_1)=
\lbrace 0\rbrace}, il suffit d'apr\`es \cite{dr1}
de prouver que \ \m{\mu(F(-2))\leq\mu(G_1)}.
Si \ \m{F=\O} \ c'est \'evident car \ \m{G_1=\O(-1)}. Sinon, on a
\ \m{\mu(G_1)-\mu(G_0)\geq 2}, et si \ \m{\mu(F(-2))>\mu(G_1)}, on a \
\m{\mu(F)-\mu(G_0)>4}, ce qui est faux car \ \m{\mu(F)-\mu(G_0)<3}.

Montrons maintenant que \ \m{\Hom(\E,G_0(1))=\lbrace 0\rbrace}. On a une
suite exacte
$$(\Hom(F,G_0(1))\ot\cx{k})\oplus(\Hom(G_0(3),G_0(1))\ot\cx{m_2}) 
\ \ \ \ \ \ \ \ \ \ \ \ \ \ \ $$
$$\ \ \ \ \ \ \ \ \ \ \ \ \ \ \ \lra
\Hom((\E,G_0(1))\lra\Ext^1(G_1(3),G_0(1))\ot\cx{m_1}.$$
On a \ \m{\Hom(F,G_0(1))=\lbrace 0\rbrace} \ car \ 
\m{\mu(F)>\mu(G_1)\geq\mu(G_0(1))}, et 
\m{\Hom(G_0(3),G_0(1))=\lbrace 0\rbrace}. Il reste \`a prouver que \
\m{\Ext^1(G_1(3),G_0(1))=\lbrace 0\rbrace}. On a
$$\Ext^1(G_1(3),G_0(1)) \ \simeq \ \Ext^1(G_0(1),G_1)^* \ = \
\lbrace 0\rbrace$$
d'apr\`es \cite{dr1} et le fait que \ \m{\mu(G_0(1))\leq\mu(G_1)}. \carre

\subsection{d\'emonstration du th\'eor\`eme A}

Soient $F$ un fibr\'e exceptionnel, $r$, \m{c_1}, \m{c_2} des entiers tels que
\m{\mu(F)-x_F<\mu<\mu(F)+x_F}, \m{\Delta<\delta(\mu)} \ et \
\m{(\mu,\Delta)\not=(\mu(F),\Delta(F))}. On peut se limiter au cas o\`u
 \m{\mu(F)-x_F<\mu\leq\mu(F)}, l'autre cas s'en d\'eduisant par dualit\'e.
On a alors
$$p \ = \ r.rg(F)(P(\mu-\mu(F))-\Delta-\Delta(F)) \ \ > \ \ 0.$$
Supposons que \ \m{\mu > \delta'(\mu)}. Alors on a \ \m{p. rg(F) < r}. En effet,
ceci \'equivaut \`a 
$$\delta(\mu)-\Delta \ < \ \q{1}{rg(F)^2}$$
(cf. la figure de l'Introduction). Il existe donc des entiers \m{r'}, \m{c'_1},
\m{c'_2}, tels que $r$, \m{c_1} et \m{c_2} soient le rang et le classes de
Chern d'une somme directe d'un fibr\'e vectoriel $\U$ de rang \m{r'} et de
classes de Chern \m{c'_1},\m{c'_2} et de \m{F\ot\cx{p}}. Le point correspondant
\`a $\U$ est situ\'e sur la conique d'\'equation
$$\Delta = P(\mu-\mu(F))-\Delta(F)$$
et on a \ \m{\Delta \geq \delta'(\mu)} \ si et seulement si ce point est situ\'e
sur le segment \m{G(F)} de la conique.

Supposons que \ \m{\Delta \geq \delta'(\mu)} \ et \ \m{r'\geq 2}. Dans ce
cas il existe d'apr\'es le lemme 3.1 un fibr\'e stable $\U$ de rang \m{r'} et de
classes de Chern \m{c'_1},\m{c'_2} tel que \ \m{\Ext^1(\U,F)=\lbrace 0\rbrace}.
Le fibr\'e 
$$\E \ = \ (F\ot\cx{p})\oplus\U$$
est prioritaire, de rang \m{r} et de classes de Chern \m{c_1},\m{c_2}. Les
fibr\'es prioritaires g\'en\'eriques sont de ce type, car les fibr\'es tels que 
$\E$ sont d\'efinis par la suite de conditions ouvertes suivante :

\medskip

\noindent (i) on a \ \m{\Ext^2(F,\E)=\lbrace 0\rbrace}.

\noindent (ii) Le morphisme canonique d'\'evaluation
$$ev : F\ot\cx{p}=F\ot\Hom(F,\E)\lra\E$$
est injectif.

\noindent (iii) Si \ \m{\U=\coker(ev)}, $\U$ est un fibr\'e stable tel que
\ \m{\Ext^1(\U,F)=\lbrace 0\rbrace}.

\medskip

Supposons maintenant que \ \m{r'=1}. Dans ce cas on doit avoir \ \m{F=\O} \ et
\m{c_2=1}. Les faisceaux de \m{M(r',c'_1,c'_2)} sont de la forme \m{\I_x} 
(id\'eal d'un point $x$ de \proj{2}). On a \ \m{\Ext^1(\I_x,\O)=\cx{}}, d'o\`u
le th\'eor\`eme A dans ce cas.

Il reste \`a traiter le cas o\`u \ \m{\Delta < \delta'(\mu)}. C'est une cons\'equence du th\'eor\`eme C, dont la d\'emonstration suit. \carre

\subsection{\D\'emonstration du th\'eor\`eme C}

Soit \m{(E,F,G)\in{\bf T}}. En consid\'erant la suite spectrale de Beilinson 
g\'en\'eralis\'ee associ\'ee \`a \m{(E,F,G)}, on voit imm\'ediatement que les
points \m{(\mu,\Delta)} de \m{\T_{(E,F,G)}} (\`a coordonn\'ees rationnelles)
sont les paires \m{(\mu(\E),\Delta(\E))}, o\`u $\E$ est de la forme
$$\E\ = (E\ot\cx{a})\oplus(F\ot\cx{b})\oplus(G\ot\cx{c}),$$
avec \ $a,b,c\geq 0$ \ non tous nuls. Le fibr\'e pr\'ec\'edent est 
prioritaire et rigide, c'est donc un fibr\'e prioritaire g\'en\'erique.

On pose comme dans le lemme 2.3, 
$$Z \ = \ \bigcup_{(E,F,G)\in{\bf T}}\T_{(E,F,G)}.$$
La partie 1- du th\'eor\`eme C est une cons\'equence imm\'ediate du \para 2.4.
Il reste donc \`a prouver que 
$$Z \ = \ \SS.$$
Soit \ \m{(\mu,\Delta)\in Z}. Alors on a \ \m{\Delta\leq\delta'(\mu)}, car
les fibr\'es prioritaires g\'en\'eriques ayant les invariants \m{\mu} et
\m{\Delta} sont rigides, comme on vient de le voir. On a donc \
\m{Z\subset\SS}.

Soit $F$ un fibr\'e exceptionnel tel que \ \m{-1<\mu(F)\leq 0}, \m{(G_n)} la
s\'erie exceptionnelle \`a gauche associ\'ee \`a $F$. On va montrer que lorsque
$n$ tend vers l'infini, le segment de conique \m{G_nF} de 
\m{T_{(G_{n-1},G_n,F)}} tend vers le segment de conique
$$\lbrace(\mu,\delta'(\mu)), \mu(F)-x_F<\mu\leq\mu(F)\rbrace.$$
On montrerait de m\^eme que si  \ \m{-1\leq\mu(F)<0}, et si \m{(H_n)} est la
s\'erie exceptionnelle \`a droite associ\'ee \`a $F$, alors lorsque
$n$ tend vers moins l'infini, le segment de conique \m{FH_n} de 
\m{T_{(F,H_n,H_{n+1})}} tend vers le segment de conique
$$\lbrace(\mu,\delta'(\mu)), \mu(F)\leq\mu<\mu(F)+x_F\rbrace.$$
D'apr\`es le lemme 2.3, ceci entraine que \ \m{\SS\subset Z}.

L'\'equation du segment de conique \m{G_nF} de \m{T_{(G_{n-1},G_n,F)}} est
$$\Delta\ = \ P(\mu-\mu(G_{n-1})-3)-\Delta(G_{n-1}).$$
On a
$$\lim_{n\rightarrow\infty}(\mu(G_{n-1})) \ = \ \mu(F)-x_F, \ \ \ 
\lim_{n\rightarrow\infty}(\Delta(G_{n-1})) \ = \ \q{1}{2}.$$
Donc le segment \m{G_nF} tend vers la courbe
$$\lbrace(\mu,\phi(\mu)), \mu(F)-x_F<\mu\leq\mu(F)\rbrace.$$
avec
$$\phi(\mu) \ = \ P(\mu-\mu(F)+x_F-3)-\q{1}{2}.$$
On v\'erifie imm\'ediatement que \ \m{\phi(\mu)=\delta'(\mu)}, ce qui ach\`eve
la d\'emonstration du th\'eor\`eme C. \carre

\end{document}